\documentclass{amsart}
\usepackage{amssymb}
\usepackage{amsmath}
\usepackage{amsfonts}
\theoremstyle{plain}
\newtheorem{thm}{Theorem}[section]
\newtheorem{prop}[thm]{Proposition}
\newtheorem{lmm}[thm]{Lemma}

\newtheorem{qws}[thm]{Question}

\theoremstyle{remark}

\def\RR{\mathbb{R}}

\def\II{\mathbb{I}}

\def\pmc#1{\setbox0=\hbox{#1}
    \kern-.1em\copy0\kern-\wd0
    \kern.1em\copy0\kern-\wd0}

\def\op{\operatorname}
\def\ov{\overline}

\def\sm{\setminus}
\begin{document}

\title[ The fundamental group $\pi_{1}(\mathbb R^3 / C)$ is uncountable] {On the fundamental group of
$\mathbb R^3$ modulo the Case-Chamberlin continuum}
\author{Katsuya Eda}
\address{School of Science and Engineering,
Waseda University, Tokyo 169-8555, Japan}
\email{eda@logic.info.waseda.ac.jp}
\author{Umed H. Karimov}
\address{Institute of Mathematics,
Academy of Sciences of Tajikistan,
Ul. Ainy $299^A$, Dushanbe 734063, Tajikistan}
\email{umed-karimov@mail.ru}
\author{Du\v san Repov\v s }
\address{Institute of Mathematics, Physics and Mechanics, and Faculty of Education,
University of Ljubljana, P.O.Box 2964,
Ljubljana 1001, Slovenia}
\email{dusan.repovs@guest.arnes.si}

\date{May 29, 2007}


\begin{abstract}It has been known for a long time that the fundamental group of the quotient of $\RR ^3$ by the
 Case-Chamberlin continuum is nontrivial. In the present paper we prove that this group is in fact, uncountable.
\end{abstract}

\maketitle

\section{Introduction}
In the 1960's, during the early days of the decomposition theory, the quotient space $X^3$ of the Euclidean 3-space
$\RR^3$ by the 
classical
Case-Chamberlin continuum $C$ (see \cite{CC}) was one of the most interesting examples. One of the most important questions was whether
$X^3$ is simply connected. It was settled -- in the negative -- by Armentrout \cite{Ar} and Shrikhande \cite{S}. However, it remained  an open problem until present day to determine how big is the fundamental group of $X^3$.
In this paper we give the solution for this problem -- namely, we show that the fundamental group 
$\pi_{1}(\RR^3 / C)$
is {\sl uncountable}.

Consider the Case-Chamberlin inverse sequence $\mathcal P$ (see \cite{CC}, \cite[p.628]{KR1}): $$P_0 \stackrel{f_0}{\longleftarrow} P_1
\stackrel{f_1}\longleftarrow P_2 \stackrel{f_2}\longleftarrow
\cdots$$
\noindent 
where $P_0 = \{ p_0\}$ is a 
singleton, 
$P_i$ is a bouquet of two
circles $S^1_{a_i} \bigvee S^1_{b_i}$, and
$p_i$ is the base point of
the bouquet $S^1_{a_i} \bigvee S^1_{b_i}$, for every $i > 0.$ 

Fix an orientation on each of the circles of the bouquet. 
Let 
$$f_i
:S^1_{a_{i+1}} \bigvee S^1_{b_{i+1}} \to S^1_{a_i} \bigvee
S^1_{b_i}$$
be a piecewise linear mapping which maps the base point
$p_{i+1}$ to the base point $p_{i}$ and
maps the natural generators
$a_{i+1}$ and $b_{i+1}$ of 
$\pi_1(S^1_{a_{i+1}} \bigvee S^1_{b_{i+1}})$ to the commutators
$[a_i, b_i]$ and $[a^2_i, b^2_i]$ of  $\pi_1(S^1_{a_i} \bigvee
S^1_{b_i})$, respectively. 

The {\sl Case-Chamberlin continuum} $C$ is then defined as the
inverse limit $\lim_{\leftarrow} \mathcal P$ 
of the Case-Chamberlin inverse sequence $\mathcal P$ (see \cite{CC}). 
Obviously, $C$ is a $1$-dimensional continuum and therefore 
it is embeddable in $\RR ^3$ (see \cite{HW}). 
It is well-known that the homotopy types of the quotient space
$\RR ^3/f(C)$ are the same for all embeddings $f$ of $C$ into $\RR ^3$
(see \cite{Borsuk}). The main result of our paper is the following theorem:

\begin{thm}\label{thm:main}
Let $C$ be the Case-Chamberlin continuum embedded in $\RR^3$. 
Then the fundamental group $\pi_{1}(\RR^3 / C)$ of the quotient
space $\RR^3/C$ is uncountable.
\end{thm}

\section{Preliminaries}

Let $G$ be a group. 
By the {\sl commutator} of the
elements
$a$ an $b$ of $G$ we mean the element $[a, b] = a^{-1}b^{-1}ab$ of
$G$.
Let $G_n$ be the {\sl lower central series} 
which is defined
inductively
(see \cite{MKS}): 
$$G_1 = G,\  G_{n+1} = [G_n, G],$$
where $[G_n, G]$ is
the group generated by the set $\{[a, b]: a\in G_n, b\in G \}$.

Obviously, $G_n \supseteq G_{n+1}$, for every $n$. By the {\sl weight}
$w(g)$ of an element $g\in G$ we mean the maximal number $n$ such that $g\in
G_n$ if such a number exists, and $\infty$ otherwise. 
 So the weight of any element of a perfect group is equal to $\infty$. We shall need the following result 
 from \cite[Ch. I, Proposition 10.2]{LS}:

\begin{prop}\label{prop:weight}:
For any free group $F$ the lower central series $F_n$ has trivial
 intersection, i.e. $\bigcap _{n=1}^\infty F_n = \{ e\}$. 
\end{prop}
That is, in any free group the weight of an element $x$ is finite if and
only if $x \neq e$.
Let $$C(f_0, f_1, f_2, \dots)$$ be the infinite mapping cylinder of
$\mathcal P $ (see e.g. \cite{K, Si}) and let
$\widetilde{\mathcal P}$ be its natural compactification by 
the Case-Chamberlin continuum $C$.
Let $\mathcal P^*$ be the quotient space of $\widetilde{\mathcal
P}$ by $C$. 

Obviously, 
$\mathcal P^*$
is homeomorphic to the
one-point compactification of an infinite 2-dimensional polyhedron
$C(f_0, f_1, f_2, \dots)$. Let 
$$C(f_k, f_{k+1}, f_{k+2}, \dots)$$
be 
the
mapping cylinder of the inverse sequence:

$$P_k \stackrel{f_k}{\longleftarrow} P_{k+1}
\stackrel{f_{k+1}}\longleftarrow P_{k+2}
\stackrel{f_{k+2}}\longleftarrow \cdots.$$

We shall denote the
corresponding one--point compactification
by 
$$C(f_k, f_{k+1}, f_{k+2}, \dots)^*.$$ 
We shall consider $C(f_k, f_{k+1}, f_{k+2},
\dots)^*$ as a subspace of $\mathcal P^*$ and we shall denote the
compactification point  by $p^*.$

We consider $P_i$, for $i\ge 0$, as a subspace of $C(f_0, f_1,
\dots)$ and we consider $C(f_k, f_{k+1}, f_{k+2}, \dots)$, for $k\ge 0$, as a subspace
of $\widetilde{\mathcal P}$. 
Obviously, $P_1$ is a strong deformation retract of $C(f_1, f_2,
\dots)$. We have the following
homomorphism 
$$\varphi _{i+1} = (f_1\cdots
f_i)_{\sharp}: \pi_1(P_{i+1})\to \pi_1(P_1)$$
which is 
a monomorphism, since it is the composition of monomorphisms
$(f_i)_{\sharp}:\pi_1(P_{i+1})\to \pi_1(P_i).$ 
Note that for a 
fixed $i$, the elements $[a_i,
b_i]$ and $[a_i^2, b_i^2]$ are free generators of a subgroup 
$(f_i)_{\sharp}(\pi_1(P_{i+1}))$ of 
$\pi_1(P_i)$ (see Exercise 12 on p.119 of \cite{MKS}).

Since $\varphi_i$ is a monomorphism, we can consider the group
$\pi_1(P_i)$ as a subgroup of $\pi_1(P_1) = F,$ where
$F$ is a free group
 on two generators $a_1$ and $b_1.$ In particular, by 
identification, we have 
$$a_2 = [a_1, b_1], a_3 = [a_2, b_2] =
[[a_1, b_1], [a_1^2, b_1^2]], \hbox{etc.}$$

Since $a_i \neq e$, the weight $w(a_i$) is
a finite number
(cf. Proposition~\ref{prop:weight} above). 
It follows by definition of $a_i$
that $w(a_i) \geq i$,
for every $i$.

Choose an increasing sequence of natural
numbers $\{n_i\}$ as follows: Let $n_0 = 1$ and $n_1 = 2$. If $n_k$ is
already defined, then let $n_{k+1}$ be any natural number such that
$n_{k+1} > w(a_{n_k})$ for $k\ge 1$. 
Then we have $a_{n_k} \notin F_{n_{k+1}}$.

Let $I_i$ be the unit segment which connects the points $p_{i+1}$ and
$p_i$ and which corresponds to the mapping cylinder of the mapping
$f_i|_{\{ p_{i+1}\}}$ of the one-point set $\{p_{i+1}\}$ to the one-point 
set $\{p_i\}$, for $i\ge 0$. 

To define a certain kind of loops we need a new notion. 
For two paths $f, g:\II\to X$ satisfying $f(1)=g(0)$, let $fg:\II\to X$
be the path defined by: 
\[ fg(s) = \left\{
\begin{array}{ll}
f(2s) & \mbox{if} \ \ 0\le s\le 1/2, \\
g(2s-1) & \mbox{if} \ \ 1/2\le s\le 1.
\end{array}\right. \]
We also let 
$$\ov{f}(s) = f(1-s) \ \hbox{for} \  0\le s\le 1.$$ 

Two paths are simply said to be {\sl homotopic}, if they are homotopic
relative to the end points.
A {\sl loop} in $X$ is a path $f:\II\to X$, satisfying $f(0)=f(1)$.
For a sequence of units and zeros 
$$\varepsilon = (\varepsilon_{1}, \varepsilon_{2},
\varepsilon_{3},\cdots), \ \ \ \varepsilon_{i} \in \{0,1\}$$ 
define a path $g_{\varepsilon}:\II \to \mathcal{P}^*$ so that the following properties hold: 

\smallskip
\noindent
(1) $g_{\varepsilon}(0) = p_1$  and $g_{\varepsilon}(1) = p*$,

\noindent
(2) $g_{\varepsilon}$ maps $[(2k-2)/(2k-1),(2k-1)/2k]$ 
homeomorphically onto $\bigcup_{i=n_{k-1}}^{n_k -1} I_i$\, starting from
$p_{n_{k-1}}$ to $p_{n_k}$ for $k\ge 1$, and 

\noindent
(3) $g_{\varepsilon}$ maps $[(2k-1)/2k, 2k/(2k+1)]$ onto
$S^1_{a_{n_k}}$ as a winding in the positive direction, if $\varepsilon
_k=1$, and $g_{\varepsilon}$ maps 
$[(2k-1)/2k, 2k/(2k+1)]$ to the point set $\{p_{n_k}\}$ constantly
otherwise, for $k\ge 1$. 

\smallskip
Let $h:\II
\to \mathcal{P}^*$ be a path from $p^*$ to $p_1$ which maps $\II$
homeomorphically onto $\bigcup_{i=1}^{\infty} I_i\cup \{ p^*\}$. Finally,
let $f_{\varepsilon} = g_{\varepsilon}h$. Then $f_{\varepsilon}$ is 
a loop with base point $p_1$ corresponding to 
$${a_{\varepsilon}} = 
{a_{n_{1}}}^{\varepsilon_{1}}{a_{n_{2}}}^{\varepsilon_{2}}
{a_{n_{3}}}^{\varepsilon_{3}}\cdots .$$

\section{Proof of Theorem~\ref{thm:main}}
For our proof of Theorem~\ref{thm:main} we shall need the following two lemmata:

\begin{lmm}\label{prop:homotopy}
Let $C$ be the Case-Chamberlin continuum embedded in $\RR^3$. 
Then the quotient space $\RR^3/C$
is homotopy equivalent to the
2-dimensional compactum $\mathcal P^*$.
\end{lmm}

\begin{proof} The proof is completely analogous to the proof of
the first assertion of Theorem 1.1 of \cite{KR2}
and therefore we shall omit it.
\end{proof}

\begin{lmm}\label{prop:loop}
Let $p_0,p_1,p^*$ be distinct points in a Hausdorff space $X$ and let
$f$ be
 a loop with base point $p_1$ such that $f^{-1}(\{ p_0\})$ is empty
 and $f^{-1}(\{ p^*\})$ is a singleton.  
If $f$ is null-homotopic, then there exists a loop 
 $f'$ in $X\setminus \{ p_0,p^*\}$ such that $f$ and $f'$ are homotopic
 in $X\setminus \{ p_0\}$. 
\end{lmm}

\begin{proof} Since $f$ is null-homotopic, we have a homotopy 
$F: \II\times \II \to X$ from $f$ to the
 constant mapping to $p_1$, i.e. 
 $$F(s,0) = f(s), F(s,1) = F(0,t) = F(1,t) = p_1 \
 \hbox{for} \ 
 s,t\in \II.$$ 
Let $\{ s_0\}$ be the singleton $f^{-1}(\{ p^*\})$. Let
$M$ be the connectedness component of $F^{-1}(\{ p^*\})$ containing
 $(s_0,0)$, and $O$  the connectedness
 component of 
$\II\times \II \setminus M$ containing $\II\times \{ 1\} $. 
Define $G: \II\times \II \to X$ by: 
\[ G(s,t) = \left\{
\begin{array}{ll}
F(s,t) & \mbox{if} \ \ (s,t)\in O, \\
p^* & \mbox{otherwise.} 
\end{array}\right. \]
Then $G$ is also a homotopy from $f$ to the constant mapping to $p_1$ and 
$G^{-1}(\{ p_0\})$ is contained in $O$.  

Consider $G^{-1}(\{ p^*,p_0\})\cap O$ and $\II\times \II \sm O$. By 
 definition of $M$, $G^{-1}(\{ p^*,p_0\})\cap O$ is compact and disjoint
 from $(\II\times \II \sm O)\cup \II \times \{ 0\}$. 
Using a polygonal
neighborhood of $(\II\times \II \sm O)\cup \II \times \{
 0\}$ whose closure is disjoint from $G^{-1}(\{ p^*,p_0\})\cap O$, 
we get a piecewise linear injective path $g:\II\to \II\times \II$
 such that 
 $$\op{Im}(G\circ g) \subseteq X\setminus \{ p_0,p^*\}, \  
 g(0)\in \{ 0\}\times \II,\, \
 \hbox{and} \
 g(1)\in \{ 1\}\times \II $$
 and 
 $\op{Im}(g)$ divides 
 $\II\times\II$ into two components,
 one of which contains $G^{-1}(\{ p_0\})$ and the other 
 contains $M \cup \II\times \{ 0\}$.
 We now see that $G\circ g$ is the desired loop $f'$.
\end{proof}

{\it Proof of\/} Theorem~\ref{thm:main}.
By Lemma~\ref{prop:homotopy}, it clearly
suffices to consider ${\pi_{1}}({\mathcal P}^*)$
instead of
$\pi_{1}(\RR^3 / C)$.
Suppose therefore, that the group ${\pi_{1}}({\mathcal P}^*)$ were at most
countable. 
We can assume that $p_1$" is the base point of the space $\mathcal P*$ and all 
of its subspaces considered below.
Since
the set of all sequences of units and zeros is uncountable, there
would then exist an uncountable set $E$, such that for every
$\varepsilon, \varepsilon'$ from $E$, the loops $f_{\varepsilon}$
and $f_{\varepsilon'}$ with the base point $p_1$ would be homotopy equivalent.
Fix a loop $f_{{\varepsilon}_{0}} \ ({\varepsilon}_{0} \in
E)$. 

Then every loop $f_{\varepsilon}\ov{f_{{\varepsilon}_{0}}}$ is
null-homotopic for every $\varepsilon \in E$. 
Since $\{ s: g_{\varepsilon}\ov{g_{{\varepsilon}_{0}}}(s) = p^*\}$ is a
singleton, we can apply Lemma~\ref{prop:loop} to
$g_{\varepsilon}\ov{g_{{\varepsilon}_{0}}}$. 
Since $f_{\varepsilon}\ov{f_{{\varepsilon}_{0}}}$ is homotopic to
$g_{\varepsilon}\ov{g_{{\varepsilon}_{0}}}$ in ${\mathcal P}^*
\setminus \{ p_0 \}$, we conclude that
$f_{\varepsilon}\ov{f_{{\varepsilon}_{0}}}$
is homotopic to a loop $f'_{\varepsilon}$ in ${\mathcal P}^* \setminus 
\{ p_0,p^*\}$, where the homotopy is in ${\mathcal P}^*\setminus P_0$. 

Since $E$ is uncountable and ${\mathcal P}^* \setminus \{ p_0, p^* \}$
is homotopy equivalent to the bouquet of two circles
$S^1_{a_1}\bigvee S^1_{b_1}$, that is, $\pi _1({\mathcal P}^* \setminus
\{ p_0, p^* \})$ is countable, there
exist distinct $\varepsilon$ and $\varepsilon'$ in $E$ such that
$f'_{\varepsilon}$ is homotopic to $f'_{\varepsilon'}$
in ${\mathcal P}^* \setminus \{p_0, p^*\}$ and hence in ${\mathcal
P}^*\setminus P_0$.
It follows that $f_{\varepsilon}\ov{f_{{\varepsilon}_{0}}}$ is
homotopic to $f_{\varepsilon'}\ov{f_{{\varepsilon}_{0}}}$ and
hence $f_{\varepsilon}$ is homotopic to $f_{\varepsilon'}$ 
in ${\mathcal P}^*\setminus P_0$.
Let $k$ be the minimal number such that $\varepsilon_{k} \neq
\varepsilon'_{k}$, say $\varepsilon_{k}=1$ and $\varepsilon'_{k}=0$. 
Let $Y_k$ be the quotient space of ${\mathcal
P}^*\sm P_0$ by the the closed subspace $C(f_{k+1}, f_{k+2},f_{k+3},
\dots)^*.$
Consider the projection $$q:\pi_{1} ({\mathcal P}^* \setminus
P_0) \to \pi_{1}(Y_{n_{k+1}})$$ and let
$[f_{\varepsilon}]$ and $[f_{\varepsilon'}]$ be the homotopy classes
containing $f_{\varepsilon}$ and $f_{\varepsilon'}$ respectively. 
Since $a_{n_{k+1}}, b_{n_{k+1}}\in F_{n_{k+1}}$, $F/F_{n_{k+1}}$ is a
quotient group of $\pi_{1}(Y_{n_{k+1}})$. 
Then, $q([f_{\varepsilon}]) = 
q(a_{n_1}^{\varepsilon_1})\cdots
q(a_{n_{k-1}}^{\varepsilon_{k-1}})q(a_{n_k})$ and 
$q([f_{\varepsilon'}]) = 
q(a_{n_1}^{\varepsilon_1})\cdots q(a_{n_{k-1}}^{\varepsilon_{k-1}})$. 
Since $a_{n_k} \notin F_{n_{k+1}}$, it follows that
$q(a_{n_k})$ is non-trivial and
hence $f_{\varepsilon}$ is not homotopic to
$f_{\varepsilon'}$ in ${\mathcal P}^*\setminus P_0$. 
This contradiction shows that our initial assumption was false and
therefore $\pi_{1}({\mathcal P}^*) \cong \pi_{1}(\RR^3 / C)$ is
indeed an uncountable group, as asserted. \qed

\begin{qws}
Let $C$ be the Case-Chamberlin continuum embedded in $\RR^3$. 
Is the first singular homology group with integer coefficients
$H_{1}(\RR^3/C;\mathbb Z)$ of the quotient
space $\RR^3/C$ also
uncountable?
\end{qws}

\section {Acknowledgements}
We were supported in part by the 
Japanese-Slovenian research grant
BI--JP/03--04/2, the Slovenian Research Agency
research program No. J1--6128--0101--04 and 
the Grant-in-Aid for Scientific research (C) of Japan 
No. 16540125. 
We thank the referee for comments and suggestions.


\begin{thebibliography}{99}

\bibitem{Ar}
S. Armentrout,
\emph{unpublished manuscript}.

\bibitem{Borsuk}
K. Borsuk,
\emph{On the homotopy type of some decomposition spaces},
Bull. Acad. Polon. Sci. S\'{e}r. Sci. Math. Astron. Phys. 
\textbf{18}
(1970), 235--239.

\bibitem{CC}
J. H. Case and R. E. Chamberlin,
\emph{Characterization of tree-like continua},
Pacif. J. Math.
\textbf{10} (1960), 73--84.

\bibitem{HW}
W. Hurewicz and H. Wallman, 
{\it Dimension Theory}, 
Princeton Univ. Press, Princeton, NJ. 1941.


\bibitem{KR1}
U. H. Karimov and D. Repov\v s,
\emph{On suspensions of noncontractible compacta of trivial shape},
{Proc. Amer. Math. Soc.}
\textbf{127}
(1999), 627--632.

\bibitem{KR2}
U. H. Karimov and D. Repov\v s,
\emph{On nonacyclicity of the quotient space of $\RR^3$ by the solenoid},
{Topol. Appl.}
\textbf{133}
(2003), 65--68.

\bibitem{K}
J. Krasinkiewicz, \emph{On a methods of constructing ANR-sets. An
application of inverse limits}, {Fund. Math.} \textbf{92}
(1976), 95--112.

\bibitem{LS} 
R. C. Lyndon and P. E. Schupp,
{\sl Combinatorial Group Theory},
Ergebnisse der Mathematik und ihrer Grenzgebiete, Band 89,
Springer-Verlag, Berlin 1977.

\bibitem{MKS}
W. Magnus, A. Karras and D. Solitar,
{\sl Combinatorial Group Theory},
Dover, New York, 1976.

\bibitem{S}
N. Shrikhande, \emph{Homotopy properties of decomposition spaces},
{Fund. Math.} \textbf{116} (1983), 119--124.

\bibitem{Si}
L. C. Siebenmann, \emph{Chapman's classification of shapes. A
proof using collapsing}, {Manuscripta Math.} \textbf{16}
(1975), 373--384.

\end{thebibliography}
\end{document}